\documentclass{article}

\usepackage{graphicx}

\usepackage{amsmath,amssymb,amsthm}

\usepackage{xcolor}
\usepackage[ddmmyyyy,hhmmss]{datetime}
\usepackage{currfile}
\usepackage[height=230mm,width=170mm]{geometry}

\newcommand {\caratheodory} {Carath\' eodory}

\newcommand {\er} {\mathbb R}
\newcommand {\en} {\mathbb N}

\newcommand {\zet} {\mathbb Z}

\newcommand{\cezero}[1] {\mathrm{C}(#1;\rn)}
\newcommand{\aceje}[1] {\mathrm{AC}^1(#1;\rn)}
\newcommand{\ceje}[1] {\mathrm{C}^1(#1;\rn)}

\newcommand {\modulo} {\mathbin{\rm mod}}

\newcommand {\dd} {\, {\rm d}}

\newcommand {\interior} {\mathop{\rm int}}
\newcommand {\boundary} {\partial}

\newcommand {\rn}{\mathbb{R}^n}

\newcommand{\mult}{\mathop{\rm mult}}

\newcommand{\derivative}[2] {\ensuremath{\frac{\dd {#1}}{\dd {#2}}}}

\newcommand{\suchthat}{\quad : \quad }

\newtheorem{thm}{Theorem} \newtheorem{prop}[thm]{Proposition}
\newtheorem{lm}[thm]{Lemma}

 \newtheorem{cor}[thm]{Corollary}
\newtheorem{df}[thm]{Definition} \newtheorem{rem}[thm]{Remark}
\newtheorem{exa}[thm]{Example}

 \newcommand{\xx}{\dot{x}}
\newcommand{\xxx}{\ddot{x}}

\begin{document}

\title{Multiple solutions of the Dirichlet problem in multidimensional billiard spaces}


\author{Grzegorz Gabor\\
\\
Faculty of Mathematics and Computer Science, \\
Nicolaus Copernicus University in Toru\'n, \\
Chopina 12/18, 87-100 Toruń,\\
 Poland\\
ggabor@mat.umk.pl\\
\\
Jan Tome\v cek\footnote{The second author was supported by the Grant IGA\_PrF\_2021\_008 "Mathematical Models" of the Internal Grant Agency of Palacky University in Olomouc.}\\
\\
Department of Mathematical Analysis\\
and  Applications of Mathematics, \\
Faculty of Science, Palack\' y University,\\
17. listopadu 12, 771 46 Olomouc,\\
 Czechia\\
jan.tomecek@upol.cz}

\maketitle

\begin{abstract}
Dirichlet problem in an $n$-dimensional billiard space is investigated. In particular, the system of ODEs $\ddot x(t) = f(t,x(t))$ together with Dirichlet boundary conditions $x(0) = A$, $x(T) = B$ in an $n$-dimensional interval $K$ with elastic impact on the boundary of $K$ is considered. The existence of multiple solutions having prescribed number of impacts with the boundary is proved. As a consequence the existence of infinitely many solutions is proved, too. The problem is solved by reformulating it into non-impulsive problem with a discontinuous right-hand side. This auxiliary problem is regularized and the Schauder Fixed Point Theorem is used.
\end{abstract}



\section{Introduction and notation}\label{introduction}

Many important physical phenomena can be described by mathematical models where suitably defined trajectories move in a set of constraints, and some impacts on the boundary of this set are observed. One of the most natural and important impact law says that the angle of incidence and angle of reflection are equal. This appears, e.g., in the case of elastic collisions. The simple situation of a free motion in a bounded subset $K$ of a Euclidean space with the above impact law is known as a {\em Birkhoff billiard}. A monograph \cite{kt89ams} is a nice introduction to the theory of billiards and its relation to equations of dynamics. Notice that the uniform movement can be described by the trivial second order equation $\ddot x(t)=0$. When some external forces depending on the state $x$ and time $t$ are present, the dynamics in $K$ changes to
\begin{equation}\label{bill1}
\ddot x(t)=f(t,x(t)) \ \mbox{ for a.e. } t\in [0,T].
\end{equation}
The same form of a right-hand side can be achieved if we consider the billiard table which is not flat but uneven (see \cite{g16jmaa} or section \ref{main results} for some examples). With analogy to standard billiards, in this more general case \eqref{bill1} we will call the constraint set $K$ a {\em billiard space}.

In the paper we are interested in the following rules of movement:
\begin{equation}\label{bill2}
\begin{cases}
\xxx(t)=f(t,x(t)), & \text{if } t\in [0,T], x(t) \in \interior K, \\
\xx(s+)=\xx(s)+I(x(s),\xx(s)), & \text{if } x(s) \in \boundary K,
\end{cases}
\end{equation}
where $K\subset\rn$ is some bounded subset which will be specified later, and $I:K\times\rn\to\rn$ is an {\em impulse map} describing the impact law.

It is worth mentioning that problem \eqref{bill2} is a special case but one of the most interesting cases of the theory of impulsive differential equations. This theory received a lot of attention during the last decades, see e.g.\! some recent monographs \cite{a10s,gho13wg,rt15atlp}.
It enables to model systems driven by differential equations in which abrupt changes in a solution $x$ take place at certain instants $t\in S:=\{t : g(t,x(t))=0\}$. If $g(t,x)=\tilde g(t)$, the impulses at fixed (a priori known) moments occur. If $g(t,x)=t-\tau(x)$, one has the so-called {\em state-dependent impulses} with barriers in the extended phase-space. Both mentioned cases do not cover the billiard dynamics \eqref{bill2} we are interested in. In fact, for \eqref{bill2} we need in the set $S$ an independent of $t$ barrier map, e.g., $g(t,x)=G(x):=d_{\partial K}(x)$.

\begin{figure}[ht]
\centering
\includegraphics[width=10cm]{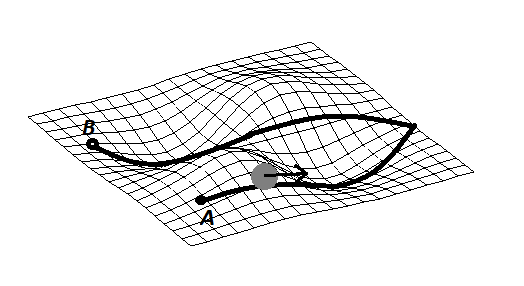}
\caption{Illustration of a two-point BVP on an uneven billiard table.}
\label{table}
\end{figure}

Among several interesting questions concerning the constraint system \eqref{bill2} there is one we would like to examine in the paper, the two-point boundary value problem with the boundary condition $x(0)=A$, $x(T)=B$, where $A,B\in K$ (comp., e.g., \cite{rt15atlp,t17ms,t17jmaa}, where the two-point boundary value problem has been examined for time- and state-dependent impulses). It is very natural when we want to hit a ball to a hole. If the billiard table is uneven or some external forces appear, both the existence and multiplicity results are nontrivial (see Figure \ref{table}).

In \cite{g16jmaa} the author obtained such results in one dimensional billiard spaces and proposed new ideas to study the problem in higher dimensions. Several open questions were stated. In response to this paper, two papers dealing with one-dimensional case were written. In \cite{abrw18anm} the authors gave the numerical treatment to the problem. In \cite{t19mmn} the author transformed a one dimensional billiard problem to the problem without impulses and proved multiplicity results by the use of the Schauder Fixed Point Theorem. It was possible because the real line could be viewed as a mosaic built of copies of an interval $K=[0,R]$, i.e., $\er=\sum_{i\in \mathbb Z} [iR,(i+1)R]$.  In the present paper we use this idea to multidimensional billiard spaces.

Namely, we examine the following two-point boundary value problem

\begin{equation}\label{billiard0}
\begin{cases}
\xxx(t)=f(t,x(t)), & \text{if } t\in [0,T], x(t) \in \interior K, \\
\xx(s+)=\xx(s)+I(x(s),\xx(s)), & \text{if } x(s) \in \boundary K, \\
x(0) = A \in \interior K,\ x(T)=B\in \interior K, &
\end{cases}
\end{equation}
where
\begin{equation}\label{K0}
K:=[\alpha_1,\beta_1]\times\ldots\times [\alpha_n,\beta_n]\subset\rn
\end{equation}
is an {\em interval in} $\rn$, and the {\em impulse map} $I(x,v)$ is given by
\[
I((x_1,\ldots,x_n),(v_1,\ldots,v_n)):= (v_1',\ldots,v_n')
\]
with $v_i'=0$ if $x_i\in (\alpha_i,\beta_i)$ and $v_i'=-2v_i$ if $x_i\in \{\alpha_i,\beta_i\}$. Therefore, the impact law can be described geometrically as the equality of angles before and after a collision with a boundary of the billiard space.

In the paper we use the following notation. By $\cezero{[a,b]}$, $\ceje{[a,b]}$, $\aceje{[a,b]}$, ${\mathrm L}^1([a,b],\rn)$ we denote spaces of all functions $x : [a,b] \to \rn$ which are continuous, have continuous derivative, have absolutely continuous derivative, are Lebesgue integrable on $[a,b] \subset \er$ respectively, for $n \in \en$.

\begin{df} \label{defi:solution}
\begin{em}
The function $x \in \cezero{[0,T]}$ is called a solution of the problem \eqref{billiard0} if and only if
\begin{itemize}
\item there exists a finite set $\{s_1,\ldots,s_p\} \subset (0,T)$, $s_1 < s_2 < \ldots < s_p$ such that $x|_{[s_i,s_{i+1}]} \in \aceje{[s_i,s_{i+1}]}$ for $i = 0,\ldots,p$, where we put $s_0 =0$ and $s_{p+1} = T$,
\item $x(t) \in K$ for each $t \in [0,T]$ and $x(t) \in \boundary K$ if and only if $t = s_i$, $i=1,\ldots,p$,
\item $x$ satisfies conditions in \eqref{billiard0}.
\end{itemize}
The number $p$ is called \emph{the number of impacts} of solution $x$ in $K$. For each $t \in [0,T]$ we define \emph{the multiplicity number of the impact at $x(t)$} as the number
\[
\mult(x,t) = \#\{k \in \{1,\ldots,n\} \suchthat x_k(t) \in \{\alpha_k,\beta_k\}\}.
\]
The number
\[
\sum_{t \in J} \mult(x,t)
\]
is called the number of impacts of the solution $x$ in $K$ up to multiplicity; we denote it by $\mult(x)$.
\end{em}
\end{df}

\begin{rem} \label{remark:number_of_impacts}
\begin{em}
The multiplicity number of the impact of the solution $x(t)$ of \eqref{billiard0} at $x(t)$ in $\er^3$ is
\begin{itemize}
    \item zero, if $x(t) \in \interior K$,
    \item one, if $x(t)$ lies inside of one of the faces of $K$,
    \item two, if $x(t)$ is at one of the edges or
    \item three, if $x(t)$ is one of the vertices.
\end{itemize}
According to Definition \ref{defi:solution}, the number of impacts $p$ of solution $x$ of \eqref{billiard0} is less or equal to $\mult(x)$. The equality holds if and only if the impacts are only at faces of $K$.\end{em}
\end{rem}

We assume that
\begin{equation}\label{ass1}
\begin{array}{l}
\alpha_i<\beta_i \mbox{ for every } i\in \{1,\ldots,n\}, \ T>0,\\
\mbox{$f:[0,T]\times K\to\rn$ is a \caratheodory\ integrably bounded map,}\end{array}
\end{equation}
i.e., $f(\cdot,x)$ is measurable for every $x\in K$, $f(t,\cdot)$ is continuous for a.e. $t\in [0,T]$, and there exists $m\in {\mathrm L}^1([0,T];\er)$ such that $||f(t,x)||\leq m(t)$ for each $(t,x)\in [0,T]\times K$. Denote  $\overline{m}:= ||m||_{L^1}$.

Notice that the study of problem \eqref{billiard0} is equivalent with the study of the following one
\begin{equation}\label{billiard}
\begin{cases}
\xxx(t)=f(t,x(t)), & \text{if } t\in [0,T], x(t) \in \interior K, \\
\xx(s+)=\xx(s)+I(x(s),\xx(s)), & \text{if } x(s) \in \boundary K, \\
x(0) = A \in \interior K,\ x(T)=B\in \interior K &
\end{cases}
\end{equation}
where
\begin{equation}\label{K}
K:=[0,c_1]\times\ldots\times [0,c_n]\subset\rn,
\end{equation}
and $I(x,v)$ is given by
\[
I((x_1,\ldots,x_n),(v_1,\ldots,v_n)):= (v_1',\ldots,v_n')
\]
with $v_i'=0$ if $x_i\in (0,c_i)$ and $v_i'=-2v_i$ if $x_i\in \{0,c_i\}$.

Indeed, if we are given, for instance, problem \eqref{billiard0}, assume that $c_i$ satisfy $c_i:=\beta_i-\alpha_i$, for $i=1,\ldots,n$. For $K$ defined in \eqref{K0} consider $\tilde K:=[0,c_1]\times\ldots\times [0,c_n]$ and $\tilde f:[0,T]\times \tilde K\to\rn$,
\[\tilde f(t,x):=f(t,x+\alpha) \ \mbox{ for a.e. } t\in [0,T] \mbox{ and } x\in \tilde K,\]
where $\alpha:=(\alpha_1,\ldots,\alpha_n)$. Consider also the points $\tilde A:=A-\alpha$ and $\tilde B:=B-\alpha$.

Notice that $f$ is integrably bounded if and only if $\tilde f$ is so, with the same integrable function $m$. Moreover, it is easy to check that $\tilde x$ is a solution to \eqref{billiard} if and only if $x(t):=\tilde x(t)+\alpha$ is a solution to \eqref{billiard0}.

Having in mind what has been said above we shall focus in what follows on problem \eqref{billiard} instead of \eqref{billiard0}.

The main result of the paper reads as follows.
\begin{thm}\label{main}
Let $A,B\in \interior K$, where $K$ is given in \eqref{K}, and $p\in \mathbb{N}$ be such that
\[p>\max_{i=1,\ldots,n}\frac{T\overline{m}}{c_i}+1.\]
Then there exist at least $2^n$ solutions of \eqref{billiard} having at least $p$ impacts and exactly $np$ impacts up to multiplicity.
\end{thm}

As a consequence (if $0\in \interior K$ and $A=B=0$) we obtain the existence and multiplicity result for the Dirichlet boundary value billiard problem with $x(0)=x(T)=0$.

The paper is organized as follows. In Section \ref{aux} we transform problem \eqref{billiard0} to the nonimpulsive second-order differential problem in $\rn$ and approximate it by a family of more regular systems. A result on the existence of a strictly monotone solution (Proposition \ref{monotone-solution}) is the main result of this section. The proof of Theorem \ref{main}, some consequences and open questions are included in Section \ref{main results}.

\section{Auxiliary problem}\label{aux}
The auxiliary problem is constructed in two steps. First, we leave the impulsive conditions at the boundary of $K$ at the cost of losing the regularity of the right-hand side of the differential equation and obtain (possibly singular) equation \eqref{aux-equation}. The possible singularity lies in the state variable, and we overcome this obstacle by constructing the sequence of regular problems \eqref{aux-equation-m}. By means of a priori bound technique and Arzel\` a--Ascoli Theorem we obtain the existence of multiple solutions of the two-point boundary value problem for singular equation \eqref{aux-equation}.

Consider, for every $i=1,\ldots,n$, the function $\theta_i:\er\to\er$,
\[\theta_i(s):=
\begin{cases}
1, & \text{if } s\modulo 2c_i\in (0,c_i),\\
-1, & \text{if } s\modulo 2c_i\in (c_i,2c_i),\\
0, & \text{if } s\modulo c_i=0,
\end{cases}\]
and
\[\Delta_i(s):=
\begin{cases}
s\modulo 2c_i, & \text{if } s\modulo 2c_i\in [0,c_i),\\
2c_i-s\modulo 2c_i, & \text{if } s\modulo 2c_i\in [c_i,2c_i).
\end{cases}\]

\begin{rem} \label{remark:delta_i}
\begin{em}
Note that for $i =1,\ldots,n$ the function $\Delta_i$ has the following properties:
\begin{itemize}
\item $\Delta_i(s) \in [0,c_i]$ for $s \in \er$,
\item $\Delta_i$ is a $2c_i$-periodic, even, continuous, piece-wise linear function.
\end{itemize}
\end{em}
\end{rem}

Consider the map $\psi:\rn\to\rn$,
\[\psi(x_1,\ldots,x_n):=(\Delta_1(x_1),\ldots,\Delta_n(x_n)).\]

We denote $f(t,x)=(f_1(t,x),\ldots,f_n(t,x))$, and extend $f$ to the map $f^*:[0,T]\times\rn\to\rn$, $f^*(t,x):=(f_1^*(t,x),\ldots,f_n^*(t,x))$, where
\[
f_i^*(t,(x_1,\ldots,x_n)):=\theta_i(x_i)f_i(t,\psi(x_1,\ldots,x_n)).
\]
It is easy to see that $f^*$ is $2c$-periodic in the sense that $f^*(t,x+2c)=f^*(t,x)$ for every $(t,x)\in [0,T]\times \rn$.

Now we are going to study the auxiliary differential equation
\begin{equation}\label{aux-equation}
\ddot z=f^*(t,z)
\end{equation}
on $[0,T]\times \rn$.

By a {\em solution} of \eqref{aux-equation} we mean a function $z\in \aceje{[0,T]}$ satisfying the equality in \eqref{aux-equation} almost everywhere. We say that a solution $z(t)=(z_1(t),\ldots,z_n(t))$ of \eqref{aux-equation} is {\em strictly monotone} if $z_i$ is a strictly monotone function for every $i=1,\ldots,n$.

\begin{prop}\label{mono-impact}
If $z$ is a strictly monotone solution of equation \eqref{aux-equation} such that $z_i(0)\modulo c_i\neq 0$ and $z_i(T)\modulo c_i\neq 0$ for every $i=1,\ldots,n$, then the function $x(t):=\psi(z(t))$ is a solution of the billiard problem \eqref{billiard} having number of impacts up to multiplicity equal to
\[
\sum_{i=1}^n \left|\left\lfloor \frac{z_i(0)}{c_i}\right\rfloor-\left\lfloor\frac{z_i(T)}{c_i}\right\rfloor\right|.\]
Moreover, the number of impacts is greater or equal to
\[
\max_{i=1,\ldots,n} \left|\left\lfloor \frac{z_i(0)}{c_i}\right\rfloor-\left\lfloor\frac{z_i(T)}{c_i}\right\rfloor\right|.
\]
\end{prop}

\begin{proof} From the definition of $x$ it follows that $x(t) \in K$ for each $t \in [0,T]$. Since $z$ is strictly increasing and continuous on $[0,T]$, there exists a finite set $M$ of instants $s \in (0,T)$ for which there exists $i \in \{1,\ldots,n\}$ such that $z_i(s)\modulo c_i = 0$. Let $J \subset (0,T) \setminus M$ be an interval. Then there exist $i_1,\ldots,i_n \in \zet$ such that for each $t \in J$ and $j =1,\ldots,n$
\[
z_j(t) \in (i_jc_j,(i_j+1)c_j)
\]
Therefore
if $i_j$ is even, then
\[
\theta_j(z_j(t)) = 1 \quad \text{and}\quad
\Delta_j(z_j(t)) = z_j(t) \modulo 2c_j = z_j(t) - i_jc_j
\]
and if $i_j$ is odd, then
$
\theta_j(z_j(t)) = -1$ and
\[\Delta_j(z_j(t)) = 2c_j - z_j(t) \modulo 2c_j = 2c_j - (z_j(t) - (i_j - 1)c_j) =
(i_j + 1)c_j - z_j(t).
\]
Let $k \in \{1,\ldots,n\}$. If $i_k$ is even, then
\[
\ddot x_k(t) = \ddot z_k(t) = f^*_k(t,z(t)) = \theta_k(z_k(t))f_k(t,x(t)) = f_k(t,x(t))
\]
and if $i_k$ is odd, then
\[
\ddot x_k(t) = -\ddot z_k(t) = -f^*_k(t,z(t)) = -\theta_k(z_k(t))f_k(t,x(t)) = f_k(t,x(t))
\]
for each $t \in J$. We have proved that $x$ is a solution of the differential equation from \eqref{billiard} on $J$.

Let $s \in M$, $j \in \{1,\ldots,n\}$. We distinguish three cases:
\begin{itemize}
\item[(a)] $z_j(s) \modulo c_j \ne 0$ or
\item[(b)] $z_j(s) \modulo 2c_j = 0$ or
\item[(c)] $z_j(s) \modulo 2c_j = c_j$.
\end{itemize}

In the first case, we get $z_j(s) \in (i_jc_j,(i_j+1)c_j)$. Therefore
\[
x_j(s) \in (0,c_j)
\]
and since the function $\Delta_j$ is also continuously differentiable at $z_j(s)$, so is $x_j$ at $s$ and therefore
\[
\dot x_j(s+) = \dot x_j(s) = \dot x_j(s-).
\]

Let us consider the case (b). Then there exists even $i_j \in \zet$ such that $z_j(s) = i_jc_j$, i.e.
\[
x_j(s) = 0.
\]
Since $z_j$ is increasing, there exists $\delta > 0$ such that
\[
\forall t \in (s-\delta,s) \quad z_j(t) \in ((i_j-1)c_j,i_jc_j)
\]
and
\[
\forall t \in (s,s+\delta) \quad z_j(t) \in (i_jc_j,(i_j+1)c_j).
\]
Then for $t \in (s-\delta,s)$ we have
\[
x_j(t) = 2c_j - (z_j(t)\modulo 2c_j) = i_jc_j - z_j(t)
\]
and for $t \in (s,s+\delta)$ we have
\[
x_j(t) = z_j(t)\modulo 2c_j =  z_j(t) - i_jc_j.
\]
Therefore
\[
\dot x_j(s+) = \dot z_j(s+) = \dot z_j(s) = \dot z_j(s-) = -\dot x_j(s-) = \dot x_j(s-) - 2\dot x_j(s-).
\]
The case (c) is similar to the case (b). Here we get
\[
x_j(s) = c_j
\]
and
\[
\dot x_j(s+) =  \dot x_j(s-) - 2\dot x_j(s-).
\]
Therefore the function $x$ satisfies the impulsive conditions from \eqref{billiard} at elements of the set $M$.
Therefore the set $M$ is the set of impulsive points of $x$ and its cardinality is the number of impacts. The estimate of number of impacts follows from the fact that if $s \in M$ there can exist more than one $i$ for which $z_i(s) \modulo c_i = 0$.
 \end{proof}

Our next task is to prove the existence of a strictly monotone solution $z$ of, possibly singular, equation \eqref{aux-equation}.

Since the right-hand side $f^*$ can have (and it really does have) discontinuity points with respect to the second variable, it is not a Carath\'eodory map. Hence, inspired by \cite{t19mmn}, we are going to regularize $f^*$. To do this, we define $\eta_n^i:[0,c_i]\to [0,1]$, by
\[
\eta_m^i(s):=
\begin{cases}
\frac{2m}{c_i}s, & \text{if } s\in [0,\frac{c_i}{2m}),\\
1, & \text{if } s\in [\frac{c_i}{2m},c_i-\frac{c_i}{2m}),\\
\frac{2m}{c_i}(c_i-s) & \text{if } s\in [c_i-\frac{c_i}{2m},c_i],
\end{cases}
\]
and consider the sequence of maps $g_m^*:[0,T]\times\rn$,
\[g_m^*(t,x):=\left(\eta_m^1(x_1\modulo c_1)f_1^*(t,x),\ldots,\eta_m^n(x_n\modulo c_n)f_n^*(t,x)\right).\]
Several observations should be made.
\begin{itemize}
    \item Each $\eta_m^i$ is continuous with a set of measure $\frac{c_i}{m}$ where $\eta_m^i$ is not equal to 1,
    \item each $\eta_m^i$ satisfies $\eta_m^i(0)=\eta_m^i(c_i)=0$,
    \item each $g_m^*$ is a Carath\'edory map which differs from $f^*$ only on a set of measure tending to 0 if $m\to\infty$,
    \item $\|g_m^*(t,x)\|\leq m(t)$ for a.e. $t\in [0,T]$ and all $x\in \rn$ and $m\in \mathbb N$.
\end{itemize}

Instead of \eqref{aux-equation} we study the sequence of equations
\begin{equation}\label{aux-equation-m}
\ddot z=g_m^*(t,z)
\end{equation}
on $[0,T]\times \rn$.

Since $g_m^*$ is more regular than $f^*$, it allows us to prove the following existence result with important estimation of velocities of solutions.

\begin{lm}\label{sol-gm}
Let us assume \eqref{ass1} and that $A=(a_1,\ldots,a_n), B=(b_1,\ldots,b_n)\in \rn$ satisfy $a_i\modulo c_i\neq 0\neq b_i\modulo c_i$.
Then, for every $m\geq 1$, there exists at least one solution $y_m\in \aceje{[0,T]}$ to \eqref{aux-equation-m} with $y_m(0)=A$, $y_m(T)=B$, and satisfying
\begin{equation}\label{bound-m}
\|y_m\|\leq K:=||A||+||B||+T\overline{m},
\end{equation}
and
\begin{equation}\label{estimate1}
\left\|{\dot y_m}(t)-\frac{B-A}{T}\right\|\leq \overline{m} \ \ \text{ for every } t\in [0,T].
\end{equation}
\end{lm}

\begin{proof} Consider the function $G:[0,T]\times [0,T]\to\er$,
\[G(t,s):=
\begin{cases}
\frac{t(s-T)}{T},& \text{for } t\leq s,\\
\frac{s(t-T)}{T},& \text{for } t>s,
\end{cases}\]
which is a standard Green function for the boundary value problem
\[\ddot y=0 \ \text{ on } [0,T], \ \ y(0)=0, \ y(T)=0.\]
We define the operator $T_m:\ceje{[0,T]}\to \ceje{[0,T]}$,
\[(T_my)(t):=\frac{t}{T}B+\frac{T-t}{T}A+\int_0^T G(t,s) g_m^*(s,y(s))\dd s.\]
Let us take the convex set
\[\Omega:=\{y\in \ceje{[0,T]} : \|y\|\leq K, \ \|\dot y\|\leq K_1\},\]
where $K$ is given in \eqref{bound-m} and $K_1=\frac{||B-A||}{T}+\overline{m}$.

Notice that, for every $t\in [0,T]$,
\[\|(T_my)(t)\|\leq \|B\|+\|A\|+\int_0^T |G(t,s)| m(s)\dd s\leq \|A\|+\|B\|+T\overline{m}=K,\]
and
\[\left\|\derivative{}{t}{(T_my)}(t)\right\|\leq \frac{\|B-A\|}{T}+\overline{m}=K_1.\]
Hence, $T_m$ maps $\ceje{[0,T]}$ into $\Omega$. Moreover, since $g_m^*$ is integrably bounded, the family $\{T_my : y\in \Omega\}$ is equibounded and uniformly continuous which implies that $T_m$ is completely continuous. Therefore we can apply the Schauder fixed point theorem to obtain a fixed point $y_m$ of $T_m$ in $\Omega$. It is easy to check that $y_m$  is a solution of \eqref{aux-equation-m} with $y_m(0)=A$ and $y_m(T)=B$. The estimate \eqref{estimate1} follows from the inequality
\begin{eqnarray*}
\left\|\dot y_m(t)-\frac{B-A}{T}\right\|&=& \left\|\derivative{}{t}{(Ty_m)}(t)-\frac{B-A}{T}\right\|\\
&\leq &
\left|\frac{\dd}{\dd t}\int_0^T G(t,s)g_m^*(s,y(s))\dd s\right|\leq \overline{m}.\end{eqnarray*}
\end{proof}

\begin{rem}
\begin{em}
There is an important meaning of the inequality in \eqref{estimate1}. It says how far velocities $\dot y_m(t)$ are from the velocity of a uniform motion.

Let us also observe that estimates \eqref{bound-m} and \eqref{estimate1} do not depend on $m$.\end{em}
\end{rem}

We are in a position to prove the existence of strictly monotone solutions of \eqref{aux-equation}.

\begin{prop}
\label{monotone-solution}
Assume that $A=(a_1,\ldots,a_n)$ and $B=(b_1,\ldots,b_n)$ satisfy $a_i\modulo c_i\neq 0\neq b_i\modulo c_i$. If $|b_i-a_i|>T\overline{m}$ for every $i=1,\ldots,n$, then there exists a monotone solution $z$ of problem \eqref{aux-equation} with $z(0)=A$ and $z(T)=B$.
\end{prop}

\begin{proof} Let, for every $m\geq 1$, $y_m\in \aceje{[0,T]}$ be a solution of \eqref{aux-equation-m} with $y_m(0)=A$, $y_m(T)=B$, and satisfying \eqref{bound-m} and \eqref{estimate1}. Of course, since $|b_i-a_i|>T\overline{m}$, we have $A\neq B$. Denote $y_m(t)=(y_m^1(t),\ldots,y_m^n(t))$. The estimate \eqref{estimate1} implies that
\begin{equation}\label{eq1-ms}
-\overline{m}+\frac{1}{T}(b_i-a_i)\leq {\dot y_m^i}(t)\leq \overline{m}+\frac{1}{T}(b_i-a_i)
\end{equation}
for every $i=1,\ldots,n$. From the assumption it follows that each function $y_m^i$ is strictly monotone (increasing or decreasing).

Note that from \eqref{bound-m} and \eqref{estimate1} it follows that $y_m$ and $\dot y_m$ are equibounded. Moreover, the estimate
\[\|\dot y_m(t)-\dot y_m(s)\|=\left\|\int_s^t \ddot y_m(\tau)\dd \tau\right\|\leq \left|\int_s^t\|g_m(\tau,y_m(\tau))\|\dd\tau\right|\leq \left|\int_s^t m(\tau)\dd \tau\right|\]
implies the uniform continuity of the sequence $(\dot y_m)$. According to Arzel\`{a}-Ascoli theorem there exists a subsequence $(y_{m_k})$ convergent to some function $z\in \ceje{[0,T]}$. Denote $z(t)=(z_1(t),\ldots,z_n(t))$. Then, passing to the limit in \eqref{eq1-ms},
\begin{equation}\label{eq1-ms-y}
-\overline{m}+\frac{1}{T}(b_i-a_i)\leq \dot{z_i}(t)\leq \overline{m}+\frac{1}{T}(b_i-a_i).
\end{equation}
Fix $i\in \{1,\ldots,n\}$. If $b_i-a_i>T\overline{m}$, then $\dot{z_i}(t)>0$ for every $t\in [0,T]$, and hence, $z_i$ is increasing. Otherwise, if $b_i-a_i<T\overline{m}$, then $\dot{z_i}(t)<0$ for every $t\in [0,T]$ which implies that $z_i$ is  decreasing. Therefore $z$ is strictly monotone.

Now, it is sufficient to check that $z$ is a solution to the differential equation \eqref{aux-equation}. To do this, let us notice that, since all $z_i$ are strictly monotone, there are finite sets $T_z^i=\{t\in (0,T) : z_i(t)\modulo c_i=0\}$. Denote $T_z:=T_z^1\cup\ldots \cup T_z^n=\{t_1,\ldots,t_p\}$ for some $p\geq 0$, where $t_1<\ldots<t_p$. Fix $j\in \{1,\ldots,p-1\}$ and take arbitrary points $s_1,s_2\in (t_j,t_{j+1})$ with $s_1\leq s_2$. Then
\begin{eqnarray*}
\dot z(s_2)-\dot z(s_1)&=&\lim_{k\to\infty}(\dot y_{m_k} (s_2)-\dot y_{m_k} (s_1))\\
&=&\lim_{k\to\infty}\int_{s_1}^{s_2}\ddot y_{m_k} (s)\dd s = \lim_{k\to\infty}\int_{s_1}^{s_2}g_{m_k}^*(s,y_{m_k}(s)))\dd s.\end{eqnarray*}
Since $(y_{m_k})$ converges to $z$ uniformly on $(t_j,t_{j+1})$ and $z$ is strictly monotone and continuous on this interval, one can find $k_0\geq 1$ such that, for each $k\geq k_0$ and for each $1\leq i\leq n$, one has
\[y_{m_k}^i(t)\modulo c_i\in \left[\frac{c_i}{2m_k},c_i-\frac{c_i}{2m_k}\right]\subset (0,c_i) \ \mbox{ for each } t\in (s_1,s_2).\]
By the very definition of $g_m^*$ we obtain
\begin{eqnarray*}
\lim_{k\to\infty}\int_{s_1}^{s_2}g_{m_k}^*(s,y_{m_k}(s)))\dd s&=& \lim_{k\to\infty}\int_{s_1}^{s_2}f^*(s,y_{m_k}(s)))\dd s\\
&=& \int_{s_1}^{s_2}f^*(s,z(s)))\dd s.\end{eqnarray*}
Hence,
\begin{equation} \label{derivative_zet}
\dot z(s_2)-\dot z(s_1)=\int_{s_1}^{s_2}f^*(s,z(s)))\dd s \quad \mbox{for each } s_1,s_2\in (t_j,t_{j+1}).
\end{equation}
Since $f^*$ is integrably bounded, $z\in \aceje{[t_j,t_{j+1}]}$ for every $1\leq j\leq p-1$, and hence $z\in \aceje{[0,T]}$.

It means that $\dot z$ is absolutely continuous which implies that it is differentiable a.e., and $\ddot z(t) = f^*(t,z(t))$ for a.e.\! $t \in (t_j,t_{j+1})$. Therefore $z$ is a solution of \eqref{aux-equation}.
Moreover, it starts at $A$ and attains $B$ at $t=T$.

 \end{proof}

\section{Main results}\label{main results}

In the previous section we have made the whole necessary preparation for the main result about the existence and multiplicity of solutions of our two-point billiard problem \eqref{billiard}. Let us finally prove Theorem \ref{main}.\\

\begin{proof}  \emph{(of Theorem \ref{main})}
Denote $U:=\{(\xi_1c_1,\ldots,\xi_nc_n) ; \xi_i\in \{-1,1\}\}$ and take any $u\in U$. The set $U$ is a set of some vertices of the mosaic on $\rn$ made of $n$-dimensional intervals generated by the interval $K$. For each $p\in \mathbb N$ we define
\[\zeta_p=(\zeta_p^1,\ldots,\zeta_p^n):=
\begin{cases}
(b_1,\ldots,b_n), & \text{if } p \text{ is } even,\\
(c_1-b_1,\ldots,c_n-b_n), & \text{if } p \text{ is } odd.
\end{cases}\]
Consider points $z_0:=A$ and $z_T=(d_1,\ldots,d_n):= pu+\zeta_p$.

Notice that, for every $i=1,\ldots,n$, we have
\[|d_i-a_i|\geq pc_i-|\zeta_p^i-a_i|\geq (p-1)c_i>T\overline{m}.\]
Hence, by Proposition \ref{monotone-solution}, there exists at least one strictly monotone solution $z$ of \eqref{aux-equation} with $z(0)=z_0$ and $z(T)=z_T$. Of course, we have $2^n$ possibilities of choice of the point $u\in U$.

It is easy to check that for any $u\in U$ one has
\[\sum_{i=1}^n \left|\left\lfloor \frac{z_i(0)}{c_i}\right\rfloor-\left\lfloor\frac{z_i(T)}{c_i}\right\rfloor\right|=np\]
and
\[
\max_{i=1,\ldots,n} \left|\left\lfloor \frac{z_i(0)}{c_i}\right\rfloor-\left\lfloor\frac{z_i(T)}{c_i}\right\rfloor\right|\geq p.
\]

At last, for every obtained solution $z$, we take $x(t):=\psi(z(t))$, see Proposition \ref{mono-impact}. From this proposition it follows that $x$ is a solution of \eqref{billiard} with exactly $np$ impacts up to multiplicity and at least $p$ impacts. Obviously, for different $z_1,z_2$ solutions $x_1=\psi\circ z_1, x_2=\psi\circ z_2$ differ, too.
 \end{proof}

The multiplicity results are given as the following corollaries.

\begin{cor}\label{cor1}
Let assumption \eqref{ass1} be satisfied. Assume that $A,B\in \interior K$, where $K$ is given in \eqref{K}. Then there exist infinitely many solutions of the impulsive two-point boundary value problem \eqref{billiard}.
\end{cor}

As we noticed at the beginning of the paper, the whole discussion is true, in particular, for the set $K$ given in \eqref{K0}
with $\alpha_i<0$ and $\beta_i>0$ (so, with $0\in \interior K$). Hence, the following multiplicity result for the Dirichlet-type problem can be obtained.

\begin{cor}
Let assumption \eqref{ass1} be satisfied for $K$ defined in \eqref{K0} with $0\in \interior K$. Then there exist infinitely many solutions of the impulsive Dirichlet boundary value problem
\begin{equation}\label{billiard1}
\begin{cases}
\xxx(t)=f(t,x(t)), & \text{if } t\in [0,T], x(t) \in \interior K, \\
\xx(s+)=\xx(s)+I(x(s),\xx(s)), & \text{if } x(s) \in \boundary K, \\
x(0) = x(T)=0. &
\end{cases}
\end{equation}
\end{cor}

\begin{exa}
\begin{em}
Consider the uneven billiard table (e.g., with some dents) which can be described as a graph of a smooth function $V$ on the rectangle $K=[\alpha_1,\beta_1]\times [\alpha_2,\beta_2]$ (see Figure \ref{table}, where the function $V(x,y)=xye^{-x^2-y^2}$ is taken on the square $[-2,2]^2$). At each point $(x,y)$ the force acting on a ball depends on $V$ at $(x,y)$. One can find the force as follows.

Denote by $Q=(0,0,-g)$ the weight of the ball, and assume, for simplicity, that the motion is frictionless.

For a point $P=(x,y,V(x,y))$ we have $Q=F_p+F_m$, where $F_p=\lambda (V_x'(x,y),V_y'(x,y),-1)$ is a pressure force (normal to the table) and $F_m$ (tangent to the table) is a component responsible for a movement. Since $\langle F_p,F_m\rangle =0$, we can check that
\[\lambda=\frac{g}{(V_x'(x,y))^2+(V_y'(x,y))^2+1}\]
and, consequently,
\[F_p=\frac{g}{(V_x'(x,y))^2+(V_y'(x,y))^2+1}(V_x'(x,y),V_y'(x,y),-1).\]
Then
\begin{eqnarray*}
F_m=Q-F_p=&\left(\frac{-g V_x'(x,y)}{(V_x'(x,y))^2+(V_y'(x,y))^2+1}, \frac{-g V_y'(x,y)}{(V_x'(x,y))^2+(V_y'(x,y))^2+1},\right. \\
&\left. -g+\frac{g}{(V_x'(x,y))^2+(V_y'(x,y))^2+1}\right).\end{eqnarray*}

Its horizontal component is then equal to
\[f(x,y):= \left(\frac{-g V_x'(x,y)}{(V_x'(x,y))^2+(V_y'(x,y))^2+1}, \frac{-g V_y'(x,y)}{(V_x'(x,y))^2+(V_y'(x,y))^2+1}\right).\]

Obviously, $f$ is integrably bounded. From Corollary \ref{cor1} it follows that for any two points $A,B\in \interior K$ there exist infinitely many solutions of the impulsive two-point boundary problem \eqref{billiard0} in $K$ with the above right-hand side $f$.
\end{em}
\end{exa}

\noindent {\bf Concluding remarks.} As we mentioned in Section \ref{introduction}, in \cite{g16jmaa} the author obtained some preliminary results on the Dirichlet problem for multidimensional billiard spaces. The technique used in \cite{g16jmaa} needed more regularity on the set $K$ (it was a smooth manifold). In contrast to that approach in the present paper we have considered a nonsmooth set $K$ and a completely different technique. There are still a lot of interesting questions for future research, i.e.,
\begin{itemize}
    \item How far (for what other billiard spaces) can one develop the technique used in the present paper?
    \item How far can one weaken regularity of the billiard space under the use of methods of \cite{g16jmaa}?
    \item For simplicity and clarity of the method used in the paper we have assumed a nice regularity of the right-hand side $f$ (see \eqref{ass1}). How far can we generalize the system?
\end{itemize}

\noindent {\bf Aknowledgement}

\noindent The authors are thankful for the anonymous referee for her/his valuable remarks.

\noindent {\em Data availability statement:} Data
sharing not applicable to this article as no datasets were generated or
analysed during the current study.

\bibliographystyle{plain}
\bibliography{all}

\end{document}